\documentclass[preprint]{elsarticle}
\usepackage[T1]{fontenc}
\usepackage[latin9]{inputenc}
\usepackage{geometry}
\usepackage{amsmath}
\usepackage{amstext}
\usepackage{amssymb}

\makeatletter
\newtheorem{thm}{Theorem}
\newtheorem{pf}{Proof}
\makeatother

\begin{document}

\title{A generalized Isserlis theorem for location mixtures\\ of Gaussian
random vectors}

\author{C.Vignat\corref{cor1}} 
\address{E.P.F.L., L.T.H.I., Lausanne, Switzerland}
\cortext[cor1]{on leave from L.S.S., Sup\'{e}lec, Orsay, France}
\ead{christophe.vignat@epfl.ch}

\begin{abstract}
In a recent paper, Michalowicz et al. provide an extension of Isserlis
theorem to the case of a Bernoulli location mixture of a Gaussian
vector. We extend here this result to the case of any location mixture
of Gaussian vector; we also provide an example of the Isserlis theorem
for a {}``scale location'' mixture of Gaussian, namely the $d-$dimensional
generalized hyperbolic distribution. 
\end{abstract}
\maketitle

\section{Introduction}

Isserlis theorem, as discovered by Isserlis \cite{Isserlis} in 1918,
allows to express the expectation of a monomial in an arbitrary number
of components of a zero mean Gaussian vector $X\in\mathbb{R}^{d}$
in terms of the entries of its covariance matrix only. Before providing
in Thm \ref{thm:thm1} the slightly generalized version of Isserlis
theorem due to Withers \cite{Withers}, we introduce the following
notations: for any set $A=\left\{ \alpha_{1},\dots,\alpha_{N}\right\} $
of integers such that $1\le\alpha_{i}\le d$ and any vector $X\in\mathbb{R}^{d}$,
we use the multi-index notation and denote
\[
X_{A}=\prod_{\alpha_{i}\in A}X_{\alpha_{i}}
\]
with the convention that for the empty set 
\[
X_{\emptyset}=1.
\]
A pairing in a set $A$ is a partition of $A$ into disjoint pairs.
We denote by $\Pi\left(A\right)$ the set of all pairings $\sigma$
in $A$: note that $\Pi\left(A\right)$ is empty if $A$ has an odd
number of elements. For a given $\sigma\in\Pi\left(A\right)$, we
denote by $A/\sigma$ the set $\left\{ i;\,\sigma=\left(i,\sigma\left(i\right)\right)\right\} ;$
finally, $\underset{A}{\sum\prod}\mathbb{E}\left(X_{i}X_{j}\right)$
denotes the sum
\[
\sum_{\sigma\in\Pi\left(A\right)}\prod_{i\in A/\sigma}\mathbb{E}\left(X_{\alpha_{i}}X_{\alpha_{\sigma\left(i\right)}}\right).
\]
In other words, for a given pairing $\sigma$ in the set $A,$ we
compute the product of all possible moments $\mathbb{E}\left(X_{i}X_{j}\right)$
where $i$ and $j$ are paired by $\sigma$; then, $\underset{A}{\sum\prod}\mathbb{E}\left(X_{i}X_{j}\right)$
denotes the sum of these products over all possible pairings in $A.$
As an example
\[
\underset{\left\{ 1,1,2,4\right\} }{\sum\prod}\mathbb{E}\left(X_{i}X_{j}\right)=\mathbb{E}\left(X_{1}^{2}\right)\mathbb{E}\left(X_{2}X_{4}\right)+2\mathbb{E}\left(X_{1}X_{2}\right)\mathbb{E}\left(X_{1}X_{4}\right).
\]
A general form of Isserlis theorem, due to Withers, is as follows.
\begin{thm}
\label{thm:thm1}If $A=\left\{ \alpha_{1},\dots,\alpha_{2N}\right\} $
is a set of integers such that $1\le\alpha_{i}\le d,\,\,\forall i\in\left[1,2N\right]$
and $X\in\mathbb{R}^{d}$ is a Gaussian vector with zero mean then
\begin{equation}
\mathbb{E}X_{A}=\underset{A}{\sum\prod}\mathbb{E}\left(X_{i}X_{j}\right)\label{eq:general-1}
\end{equation}

Moreover, if $A=\left\{ \alpha_{1},\dots,\alpha_{2N+1}\right\} $
then, under the same assumptions,
\[
\mathbb{E}X_{A}=0.
\]
\end{thm}
For example, choosing $\alpha_{i}=i,\,\,1\le i\le4$ yields the well-known
identity
\[
\mathbb{E}\left(X_{1}X_{2}X_{3}X_{4}\right)=\mathbb{E}\left(X_{1}X_{2}\right)\mathbb{E}\left(X_{3}X_{4}\right)+\mathbb{E}\left(X_{1}X_{3}\right)\mathbb{E}\left(X_{2}X_{4}\right)+\mathbb{E}\left(X_{1}X_{4}\right)\mathbb{E}\left(X_{2}X_{3}\right).
\]
However, indices $\alpha_{i}$ need not be distinct: for example,
choosing $\alpha_{i}=1,\,\,1\le i\le4$ yields
\[
\mathbb{E}X_{1}^{4}=3\mathbb{E}X_{1}^{2}.
\]

Several extensions of this result have been provided recently: in
\cite{Withers}, Withers extends Isserlis theorem to the case of noncentral
Gaussian vectors and relates the result with multivariate Hermite
polynomials; in \cite{Vignat}, a general formula for Gaussian scale
mixtures, and more generally for elliptically distributed vectors
is derived; it is applied to the computation of moments of the uniform
distribution on the sphere. In \cite{Repetowicz}, Isserlis theorem
is extended to the computation of the moments of linear combinations
of independent Student-t vectors. In \cite{Grigelionis}, Isserlis
theorem is extended to Gaussian matrix mixtures, i.e. random vectors
of the form 
\[
X=AN
\]
where $N$ is a standard Gaussian vector in $\mathbb{R}^{d}$ and
$A$ is a $\left(d\times d\right)$ random matrix. Let us also mention
the reference \cite{Kan} where the author tackles the computational
complexity of formula (\ref{eq:general-1}), using Magnus lemma to
replace a product of $n$ variables by sums of polynomials of degree
$n$ in these variables.

Recently, Michalowicz et al. \cite{Michalowicz} addressed the case
of Gaussian location mixtures: they provided an extension of Isserlis
theorem to the case of a random vector $X\in\mathbb{R}^{d}$ with
probability density
\begin{equation}
f_{X}\left(x\right)=\frac{1}{2}\phi_{R}\left(x+\mu\right)+\frac{1}{2}\phi_{R}\left(x-\mu\right)\label{eq:Olson4.1}
\end{equation}
where 
\[
\phi_{R}\left(x\right)=\frac{1}{\vert2\pi R\vert^{\frac{1}{2}}}\exp\left(-\frac{1}{2}x^{t}R^{-1}x\right)
\]
is the $d-$variate Gaussian density with zero mean and covariance
matrix $R.$ 

In the following, we give a new and simple proof of the result by
Michalowicz et al., adopting a formalism that allows us to extend
their results to the general case of an arbitrary Gaussian location
mixture. We also provide an extension of these results to the case
of a scale-location mixture of Gaussian.

\section{Extensions of the result by Michalowicz et al.}

A key observation is that the random vector $X$ with density (\ref{eq:Olson4.1})
reads
\begin{equation}
X=\epsilon\mu+\zeta\label{eq:randomized}
\end{equation}
where $\epsilon$ is a Bernoulli random variable ($\Pr\left\{ \epsilon=-1\right\} =\Pr\left\{ \epsilon=1\right\} =\frac{1}{2}$),
$\zeta\in\mathbb{R}^{d}$ is a zero mean Gaussian vector and equality
is in the sense of distributions. This stochastic representation allows
to prove easily a generalized version of the main result of \cite{Michalowicz},
namely
\begin{thm}
\label{thm:2} If $X\in\mathbb{R}^{d}$ is distributed according to
(\ref{eq:Olson4.1}) and $A=\left\{ \alpha_{1},\dots,\alpha_{2N}\right\} $
with $1\le\alpha_{i}\le d$ then
\[
\mathbb{E}X_{A}=\sum_{k=0}^{N}\sum_{\substack{S\subset A\\
\vert S\vert=2k
}
}\mu_{S}\underset{A\backslash S}{\sum\prod}\mathbb{E}\left(\zeta_{i}\zeta_{j}\right).
\]
If $A=\left\{ \alpha_{1},\dots,\alpha_{2N+1}\right\} $ then $\mathbb{E}X_{A}=0.$
\end{thm}
The simplified proof we propose is as follows: by (\ref{eq:randomized}),
\[
\mathbb{E}X_{A}=\mathbb{E}\left(\epsilon\mu+N\right)_{A}
\]
and since the product $\left(a+b\right)_{A}$ can be expanded as
\[
\left(a+b\right)_{A}=\sum_{k=0}^{2N}\sum_{\substack{S\subset A\\
\vert S\vert=k
}
}a_{S}b_{A\backslash S},
\]
we deduce that 
\[
\mathbb{E}X_{A}=\sum_{k=0}^{2N}\sum_{\substack{S\subset A\\
\vert S\vert=k
}
}\mu_{S}\mathbb{E}\left(\epsilon^{k}\right)\mathbb{E}\left(\zeta_{A\backslash S}\right).
\]
By Isserlis theorem, the expectation of the product of an odd number
of centered Gaussian random variables $\zeta_{i}$ is equal to zero
so that this expression simplifies to
\[
\sum_{k=0}^{N}\sum_{\substack{S\subset A\\
\vert S\vert=2k
}
}\mu_{S}\mathbb{E}\left(\epsilon^{2k}\right)\mathbb{E}\left(\zeta_{A\backslash S}\right).
\]
Since $\epsilon$ is Bernoulli distributed, all its even moments are
equal to 1; moreover, since $\zeta$ has zero mean, by Isserlis theorem,
$E\left(\zeta_{A\backslash S}\right)=\underset{A\backslash S}{\sum\prod}E\left(\zeta_{i}\zeta_{j}\right)$
and we obtain
\[
\mathbb{E}X_{A}=\sum_{k=0}^{N}\sum_{\substack{S\subset A\\
\vert S\vert=2k
}
}\mu_{S}\underset{A\backslash S}{\sum\prod}\mathbb{E}\left(\zeta_{i}\zeta_{j}\right)
\]
which is the desired result. The case where $A$ has an odd number
of elements is equally simple.

We note that Theorem \ref{thm:2} can be also easily deduced using
generating functions as done in \cite[Theorem 1.1]{Withers} who proves
a version of Wick's theorem for a Gaussian vector with mean $\mu\ne0$
: choosing a Bernoulli randomized version of this mean as in (\ref{eq:randomized})
yields the result.

\section{The general case of Gaussian location mixture}

With the useful representation (\ref{eq:randomized}), we can generalize
the preceding result to any kind of location mixture of Gaussian:
namely, we consider a random vector $X\in\mathbb{R}^{d}$ that reads
\begin{equation}
X=\mu+\zeta\label{eq:general}
\end{equation}
where $\zeta$ is a zero-mean Gaussian vector in $\mathbb{R}^{d}$,
independent of the random vector $\mu\in\mathbb{R}^{d}$ with probability
distribution $F_{\mu}$; note that the vector $\mu$ may be discrete
- taking values $\mu_{i}$ with probabilities $p_{i}$ - or not, but
we don't need to assume the existence of a density $f_{\mu}.$ In
the discrete case, the density of $X$ reads
\[
f_{X}\left(x\right)=\sum_{i=0}^{+\infty}p_{i}\phi_{R}\left(x-\mu_{i}\right);
\]
and in the most general case,
\[
f_{X}\left(x\right)=\int_{\mathbb{R}^{d}}\phi_{R}\left(x-\mu\right)dF_{\mu}\left(\mu\right).
\]

We now state our main theorem.
\begin{thm}
Assume that $X\in\mathbb{R}^{d}$ follows model (\ref{eq:general})
and that all the first-order moments $m_{k}=\mathbb{E}\mu_{k}$ of
$\mu$ exist. Then if $A=\left\{ \alpha_{1},\dots,\alpha_{2N+\epsilon}\right\} ,$
with $\epsilon\in\left\{ 0,1\right\} ,$ 
\begin{equation}
\mathbb{E}X_{A}=\sum_{k=0}^{N}\sum_{\substack{S\subset A\\
\vert S\vert=2k+\epsilon
}
}\mathbb{E}\left(\mu_{S}\right)\underset{A\backslash S}{\sum\prod}\mathbb{E}\left(\zeta_{i}\zeta_{j}\right)\label{eq:EXA1}
\end{equation}

\end{thm}
We remark that if all elements $\alpha_{i}$ of $A$ are different
and if the vector $\mu$ has independent components, this expression
can be further simplified to 
\begin{equation}
\mathbb{E}X_{A}=\sum_{k=0}^{N}\sum_{\substack{S\subset A\\
\vert S\vert=2k+\epsilon
}
}\left(\mathbb{E}\mu\right)_{S}\underset{A\backslash S}{\sum\prod}\mathbb{E}\left(\zeta_{i}\zeta_{j}\right),\label{eq:EXA2}
\end{equation}
noting the difference between $\mathbb{E}\left(\mu_{S}\right)=\mathbb{E}\prod_{\alpha_{i}\in S}\mu_{\alpha_{i}}$
in (\ref{eq:EXA1}) and $\left(\mathbb{E}\mu\right)_{S}=\prod_{\alpha_{i}\in S}\mathbb{E}\mu_{\alpha_{i}}$
in (\ref{eq:EXA2}). 

The proof is as follows.
\begin{pf}
From (\ref{eq:general}), we deduce
\[
\mathbb{E}X_{A}=\sum_{k=0}^{2N+\epsilon}\sum_{\substack{S\subset A\\
\vert S\vert=k
}
}\mathbb{E}\left(\mu_{S}\right)\mathbb{E}\left(\zeta_{A\backslash S}\right).
\]
Since the cardinality of $A\backslash S$ is $2N+\epsilon-k,$ $\mathbb{E}\left(\zeta_{A\backslash S}\right)=0$
unless $\vert S\vert=k$ has the same parity as $\epsilon$, in which
case it is equal to $\underset{A\backslash S}{\sum\prod}\mathbb{E}\left(\zeta_{i}\zeta_{j}\right)$,
hence formula (\ref{eq:EXA1}). Formula (\ref{eq:EXA2}) is easily
deduced from formula (\ref{eq:EXA1}) assuming that the components
of $\mu$ are independent and that all elements of $A$ are distinct.
\end{pf}
We now provide a further generalization of Isserlis theorem by considering
a Gaussian vector with both random scale and location parameters.

\section{A Normal variance-mean mixture application}

The generalized $d-$dimensional hyperbolic distribution was introduced
by Barndorff-Nielsen in 1978 \cite{Barndorff-Nielsen}. It is the
distribution of a random vector that reads
\begin{equation}
X=\text{\ensuremath{\mu}}+\sigma^{2}\Delta\beta+\sigma\Delta^{\frac{1}{2}}\zeta\label{eq:X hyperbolic}
\end{equation}
where $\mu$ and $\beta$ are two deterministic vectors in $\mathbb{R}^{d}$,
$\Delta$ is a deterministic $\left(d\times d\right)$ matrix with
$\vert\Delta\vert=1,$ $\zeta$ is a standard Gaussian vector in $\mathbb{R}^{d}$
and $\sigma^{2}$ is a scalar random variable that follows the Generalized
Inverse Gaussian $GIG\left(\psi,\chi,\lambda\right)$ distribution
\begin{equation}
f_{\psi,\chi,\lambda}\left(x\right)=\frac{\left(\frac{\psi}{\chi}\right)^{\frac{\lambda}{2}}}{2K_{\lambda}\left(\sqrt{\psi\chi}\right)}x^{\lambda-1}\exp\left(-\frac{1}{2}\left(\chi x^{-1}+\psi x\right)\right),\,\, x>0\label{eq:GIGintegral}
\end{equation}
with parameters $\psi>0,\,\,\chi>0$ and $\lambda\in\mathbb{R}.$
We note that in (\ref{eq:X hyperbolic}), the GIG random variable
$\sigma^{2}$ appears both as a scale and location parameter of the
Gaussian vector, hence the {}``normal variance-mean mixture'' name.
From the stochastic representation (\ref{eq:X hyperbolic}), we derive
a version of the Isserlis theorem as follows.
\begin{thm}
If $X\in\mathbb{R}^{d}$ is a generalized hyperbolic vector as in
(\ref{eq:X hyperbolic}) and $A=\left\{ \alpha_{1},\dots,\alpha_{2N+\epsilon}\right\} $
with $\epsilon\in\left\{ 0,1\right\} $ then
\[
\mathbb{E}X_{A}=\sum_{\begin{array}{c}
0\le l\le N\\
0\le p\le2l+\epsilon
\end{array}}\sum_{\substack{T\subset S\subset A\\
\vert T\vert=p,\vert S\vert=2l+\epsilon
}
}\mu_{T}\gamma_{S\backprime T}m_{N+l-p+\epsilon}\underset{A\backslash S}{\sum\prod}\mathbb{E}Z_{i}Z_{j}
\]
where $\gamma=\Delta\beta,$ where $Z$ is a centered Gaussian vector
with covariance matrix $\Delta$ and
\[
m_{l}=\mathbb{E}\sigma^{2l}=\left(\frac{\psi}{\chi}\right)^{-\frac{l}{2}}\frac{K_{\lambda+l}\left(\sqrt{\psi\chi}\right)}{K_{\lambda}\left(\sqrt{\psi\chi}\right)}.
\]
\end{thm}
\begin{pf}
Assuming first $\epsilon=0,$ we have
\[
\mathbb{E}X_{A}=\mathbb{E}\sum_{l=0}^{N}\sum_{\substack{S\subset A\\
\vert S\vert=2l
}
}\left(\mu+\sigma^{2}\gamma\right)_{S}\left(\sigma Z\right)_{A\backslash S}=\sum_{l=0}^{N}\sum_{\substack{S\subset A\\
\vert S\vert=2l
}
}\mathbb{E}\left(\sigma^{2N-2l}\left(\mu+\sigma^{2}\gamma\right)_{S}\right)\mathbb{E}Z_{A\backslash S}
\]
with 
\[
\mathbb{E}Z_{A\backslash S}=\underset{A\backslash S}{\sum\prod}\mathbb{E}Z_{i}Z_{j}
\]
and
\[
\mathbb{E}\sigma^{2N-2l}\left(\mu+\sigma^{2}\gamma\right)_{S}=\sum_{p=0}^{2l}\sum_{\substack{T\subset S\\
\vert T\vert=p
}
}\mu_{T}\gamma_{S\backslash T}\mathbb{E}\left(\sigma^{2}\right)^{N+l-p}.
\]
The moment of order $l$ of the GIG random variable $\sigma^{2}$
can be easily computed from (\ref{eq:GIGintegral}) as
\[
m_{l}=\left(\frac{\psi}{\chi}\right)^{-\frac{l}{2}}\frac{K_{\lambda+l}\left(\sqrt{\psi\chi}\right)}{K_{\lambda}\left(\sqrt{\psi\chi}\right)},
\]
hence the result. The case $\epsilon=1$ follows the same steps.\end{pf}

\end{document}